\documentclass{amsart}

\usepackage{comment, hyperref, amsmath, amssymb}
\begin{document}

\newtheorem{theorem}{Theorem}[section]
\newtheorem{prop}[theorem]{Proposition}
\newtheorem{lemma}[theorem]{Lemma}
\newtheorem{cor}[theorem]{Corollary}
\newtheorem{conj}[theorem]{Conjecture}

\theoremstyle{definition}
\newtheorem{definition}[theorem]{Definition}
\newtheorem{rmk}[theorem]{Remark}
\newtheorem{eg}[theorem]{Example}
\newtheorem{qn}[theorem]{Question}
\newtheorem{defn}[theorem]{Definition}

\numberwithin{equation}{section}

\newcommand\Z{{\mathbb{Z}}}
\newcommand\al{{\alpha}}

\newcommand{\C}{{\mathbb C}}
\newcommand\R{{\mathbb R}}

\newcommand{\bbar}{\overline}
\newcommand{\til}{\widetilde}
\newcommand{\Ra}{\longrightarrow}

\title[A splitting theorem for good complexifications]{A splitting theorem for good 
complexifications}

\author[I. Biswas]{Indranil Biswas}
\address{School of Mathematics, Tata Institute of Fundamental
Research, Homi Bhabha Road, Bombay 400005, India}
\email{indranil@math.tifr.res.in}

\author[M. Mj]{Mahan Mj}
\address{RKM Vivekananda University, Belur Math, WB 711202, India}
\email{mahan.mj@gmail.com, mahan@rkmvu.ac.in}

\author[A. J. Parameswaran]{A. J. Parameswaran}
\address{School of Mathematics, Tata Institute of Fundamental
Research, Homi Bhabha Road, Bombay 400005, India}
\email{param@math.tifr.res.in}



\date{}

\thanks{The first and second authors acknowledge the support of their respective J. C. Bose Fellowships.}

\begin{abstract}
The purpose of this paper is to produce restrictions on fundamental groups of manifolds admitting
good complexifications by proving the following Cheeger-Gromoll type splitting theorem:\\
Any closed manifold $M$ admitting a good complexification has a finite-sheeted regular
covering $M_1$ such that $M_1$ admits a fiber bundle structure with base $(S^1)^k$ and
fiber $N$ that admits a good complexification and also has zero virtual first
Betti number. We give several applications to manifolds of dimension at most 5.
\end{abstract}

\maketitle

\tableofcontents

\section{Introduction}\label{se1}

A good complexification of a closed smooth manifold $M$ is defined to be a smooth affine algebraic variety
$U$ over the real numbers such that $M$ is diffeomorphic to $U (\R)$ and the inclusion $U (\R) \,\Ra\, U (\C)$ is a
homotopy equivalence \cite{tot}, \cite{kul}. A good complexification comes naturally
equipped with a natural antiholomorphic involution $A$ on $U (\C)$ whose fixed point set is precisely the set
of real points $U (\R)$. Kulkarni \cite{kul} and Totaro \cite{tot} investigate the topology of good
complexifications using characteristic classes and Euler characteristic. In this paper
we prove a Cheeger-Gromoll type splitting theorem
and initiate a systematic study of fundamental groups of good
complexifications.

In \cite[p. 69, 2nd para]{tot}, Totaro asks the following question:

\begin{qn}\label{tqn}
If a closed smooth manifold $M$ admits a good
complexification, does $M$ also admit a metric of non-negative
curvature?
\end{qn}

The following Cheeger--Gromoll type splitting theorem is the main result proved here.

\begin{theorem}[{See Theorem \ref{albbdl}}]\label{mt} Let $M $ be a closed manifold admitting a good complexification.
Then $M$ has a finite-sheeted regular covering $M_1$ satisfying the following:
\begin{enumerate}
\item $M_1$ admits a fiber bundle structure with fiber $N$ and base $(S^1)^d$. Here $d$ denotes the (real) Albanese dimension of $M_1$.

\item The first virtual Betti number $vb_1(N)\,=\, 0$,

\item $N$ admits a good complexification.
\end{enumerate}
\end{theorem}

Let us spell out the analogy with the Cheeger-Gromoll splitting theorem  \cite[Theorem 69 , page 288]{pet}: \\
A closed manifold of non-negative curvature has a finite sheeted cover of the form $N \times (S^1)^d$, where $N$ is simply connected.
Theorem \ref{mt} similarly furnishes a fibering over a torus $(S^1)^d$ with fiber $N$ having   $vb_1(N)\,=\, 0$. 

\smallskip

Gromov \cite{gromov-cmh} proves that if a closed smooth manifold $M$ of dimension $n$ admits a metric of non-negative
curvature, then there is an upper bound, that depends only on $n$, on the sum of the Betti numbers of $M$. He further
conjectures, \cite[Section 7]{gromov-qns},
that $b_i(M)\,\leq\, b_i((S^1)^n)$. Theorem \ref{gcg0} below furnishes positive evidence towards a combination of
Question \ref{tqn} with this conjecture of Gromov by giving an affirmative answer for the first Betti number of manifolds
admitting good complexifications.
 
We shall say that a finitely presented group $G$ is a {\bf good complexification group} if $G$ can be realized as the
fundamental group of a closed smooth manifold admitting a good complexification (see also \cite{abckt}). 
We deduce from Theorem \ref{mt} the following
critical restriction on good complexification groups (See Theorem \ref{gcg}):

\begin{theorem} \label{gcg0}
 Let $G$ be a good complexification group. Then there exists a finite index
subgroup $G_1$ of $G$ such that two following statements hold:
\begin{enumerate}
 \item There is an exact sequence: $$1\,\longrightarrow\, H\,\longrightarrow\,
G_1 \,\longrightarrow\, \Z^k\,\longrightarrow\, 1,$$ where $k$ can be zero.

\item The above $H$ is a finitely presented good complexification group with $vb_1(H)\,=\,0$, where
$vb_1(H)$ denotes the virtual first Betti number of $H$.
\end{enumerate}
\end{theorem}

Recall that for a group $H$, the {\bf virtual first Betti number} $vb_1(H) $ is the supremum of first Betti numbers
$b_1(H_1)$ as $H_1$ runs over finite index subgroups of $H$.

The following classes of groups are then ruled out as good complexification groups (See
Section \ref{gcgrest}, especially Corollary \ref{not}):

\begin{enumerate}
 \item Groups with infinite $vb_1$, in particular large groups.
\item Hyperbolic CAT(0) cubulated groups.
\item Solvable groups that are not virtually abelian.
\item 2- and 3-manifold groups that are not virtually abelian.
\item any group admitting a surjection onto any of the above.
\end{enumerate}

Question \ref{tqn} has an affirmative answer for 2-manifolds; this is probably classical
but follows also from \cite{kul, tot}. An affirmative answer to Question \ref{tqn} for $3$-manifolds
was given in \cite{bm-imrn}. As a consequence of the above restrictions,
Question \ref{tqn} has an affirmative answer for 2 and 3-manifold (see Sections
\ref{2m} and \ref{3m}). Thus a new self-contained proof of the main Theorem
of \cite{bm-imrn} on good complexifications is obtained.
In Section \ref{ld} we give a number of applications to low-dimensional manifolds.

\begin{theorem} \label{ld0} \mbox{}

\begin{enumerate}
 \item Question \ref{tqn} has an affirmative answer for 2-manifolds.

\item Question \ref{tqn} has an affirmative answer for 3-manifolds \cite{bm-imrn}.

\item Let $M$ be a closed simply connected 4-manifold admitting a symplectic good complexification. Then $M$ admits
a metric of non-negative curvature.

\item Let $M$ be a closed 4-manifold admitting a good complexification. Further suppose that
$\pi_1(M)$ is infinite, torsion-free abelian. Then $\pi_1(M)$ is isomorphic to $\Z^d$, where
$d\,=\,1, 2$ or $4$. Moreover, the manifold $M$ admits a finite-sheeted cover with a
metric of non-negative curvature (i.e., Question \ref{tqn} has an affirmative answer
up to finite-sheeted covering).

\item Let $M$ be a closed 5-manifold admitting a good complexification. Further suppose that $\pi_1(M)$ is infinite, torsion-free
abelian. Then $\pi_1(M)$ is isomorphic to 
 $\Z^d$, where $d\,=\,1, 2, 3$ or $5$. Further, if $d \,=\, 2, 3$ or $5$, then $M$ admits
a finite-sheeted cover $M_1$ homeomorphic to $S^2 \times T^3$ or $S^3 \times T^2$ or $T^5$.
In particular, $M_1$ admits a metric of non-negative curvature.
\end{enumerate}
\end{theorem}

\medskip

For an algebraic variety $X$, we shall denote the real (respectively, complex) points by
$X_\R$ (respectively, $X_\C$).
\begin{defn}
We say that an algebraic map $f\,:\, X \,\longrightarrow\, Y$ between good complexifications is
a {\bf good complexification map} if there is a commutative diagram
\begin{equation}
\begin{matrix}
X_\R && \subset && X_\C
\\
 ~\Big\downarrow f_\R && &&
~\Big\downarrow f_\C\\
Y_\R && \subset && Y_\C
\end{matrix}
\end{equation}
such that $f_\C$ is equivariant with respect to the antiholomorphic
involutions $A_X , A_Y$, meaning $f_\C \circ A_X \,=\, A_Y \circ f_\C$.
\end{defn}

\section{Fibrations and the Albanese}\label{sectech}

We need to first introduce some terminology following \cite{tib}.
Let $f\,:\, X \,\longrightarrow\, \C^\ast$ be an algebraic map of quasiprojective
varieties. This $f $ is said to be topologically trivial at $z \,\in\, \C^\ast$ if there
is a neighborhood $D$ of $z$ such that the restriction of $f$ from $ f^{-1}(D)$ to $D$
is a topologically trivial fibration. If $z$ does not satisfy this property, then
we say that $z$ is an {\bf atypical value} and that $ f^{-1}(z)$ is
an {\bf atypical fiber}. Clearly, critical values (respectively, critical
fibers) are atypical values (respectively, atypical fibers). The book \cite{tib}
gives examples of atypical values that are not critical.

\subsection{Fibration over $\C^\ast$}\label{fibncstar}

The punctured plane $\C^\ast\,=\, \C\setminus\{0\}$ is the standard good complexification
of the unit circle $S^1$ with the antiholomorphic 
involution given by $z \,\longrightarrow\, \frac{1}{\overline{z}}$. In this section we study good complexification
maps $f_\C\,:\, X \,\longrightarrow\, \C^\ast$ restricting to 
$f_\R\,:=\, f_\C\vert_{X_\R}\, :\, X_\R \,\longrightarrow\, S^1$ and prove that such maps have to be topologically trivial fibrations. 

\begin{theorem}\label{cstar}
Let $X$ be a good complexification of $M$, and let $\C^\ast$ be the standard good
complexification of $S^1$. Let $f\,:\, X\,\longrightarrow\, \C^\ast$ be a good
complexification map between the good complexifications of $M$ and $S^1$. The
corresponding map between real loci
$$f_\R\, :=\, f\vert_{M}\,:\, M \, \longrightarrow\, S^1$$
is a smooth fiber bundle. If $N$ denotes the fiber of $f_\R$, then $N$ admits
a good complexification.
\end{theorem}

\begin{proof}
The proof will be divided into a number of steps. We fix some notation first. The
map $f$ will also be denoted by $f_\C$. The restriction of $f$ to $f^{-1}(S^1)
\,\subset\, X_\C$ will be denoted as $\widehat f$.

In what follows a {\it small perturbation} will refer to a small perturbation in the 
real algebraic category. \\

\noindent {\bf Step 1: $f_\R$ has no critical points after small perturbation:}\\
First, note that since $b_1(M) \,= \,b_1(X) \,\neq\, 0$, it follows from \cite{kul,tot} that
the Euler characteristic $\chi(M)\,=\, \chi (X)\,=\,0$. Next, assume for the time being, that
$f$ can be perturbed to a function $g$ with isolated quadratic (i.e., Morse-type)
singularities in a small tubular neighborhood $X_M$ of $M=X(\R)$ such that the image of $g$ is an annular neighborhood
$A_\C$ of $S^1$. Let $F$ denote the general fiber of $g$.

Then the Euler characteristics of $F$, $X_M$ and $A_\C$ are related by the following formula due to Suzuki \cite{Su}
(Suzuki emphasizes dimension 2, though the formula works in the general context of Stein manifolds and Stein morphisms): $$
\chi (X_M) \,=\, \chi (F)\cdot\chi (A_\C) + \sum_{x\in A_\C} (\chi(F_x) - \chi(F))\, ,$$ where $F_x$ denotes the fiber of $g$
over $x$. Let $c_x$ be  the
defect in the Euler characteristic given by $(\chi(F_x) - \chi(F))$. Since $g$ is defined in
the neighborhood of $M$, the only exceptional fibers of $g$ are critical (i.e. atypical non-critical fibers
do not exist). Let $c_i$ denote the finitely many non-zero defects.  Each isolated quadratic critical
point contributes a cell of a fixed real dimension $m$ (one more than the complex dimension of the fiber).
Hence all the $c_i$'s are either $+1$ or $-1$ depending on the parity of $m$. Since
$\chi(X)\,=\,0\, =\, \chi(\C^\ast)$, it
follows from the Euler characteristic formula above that $g$ has no critical points.

It remains to show that $f$ can indeed be perturbed to a function $g$
with isolated quadratic (i.e., Morse-type) singularities locally. Since $f$ is algebraic onto $\C^\ast$, it has finitely many
critical and atypical values. Restricted to a small tubular neighborhood of $M$, $f$ has only critical points.
Let $A_\C $ denote a small annular neighborhood of $S^1$ containing the corresponding critical values.
 Any small perturbation
of the annulus $A_\C $ is again an annulus (unlike $\C^\ast$, which can be perturbed slightly
to $\C$). We note that $f$ is an affine
and hence a Stein morphism, and $A_\C $ is Stein. $X_M$ is also Stein. 
 Since $f$ is a Stein morphism between Stein manifolds with only critical points, it can be perturbed
slightly to a Stein morphism with isolated quadratic singularities, with all singular  values in 
$A_\C$.

\noindent {\bf Step 2: Existence of critical points of $f_\C$ is stable under 
perturbation; hence $ f_\R$ has no critical points:}\\ Further, if $f_\C$ has critical points inside $X_M$, then any small
perturbation of $f_\C$ must also have critical points. This is clear for isolated 
critical points. In general, the set of critical points $Z(f_\C)$ is a variety and 
taking a local (algebraic) section $\sigma$ of $f_\C$ with the image of $\sigma$
lying in a neighborhood of $Z(f_\C)$, we can 
ensure that the intersection $Z(f_\C)\bigcap \sigma$ is an isolated point. Restricting 
$f_\C$ to $Z(f_\C)\bigcap\sigma$ we observe that a small perturbation of $f_\C$ must 
have critical points in $Z(f_\C)\bigcap\sigma$. 

Since a small perturbation $g$ of $f$ has no critical points inside $X_M$, it follows that $f=f_\C$ has no critical points
inside $X_M$.
If $f_\R$ has critical points, so does $f_\C$ inside $X_M$. It follows that $f_\R$ has no critical points.

\noindent {\bf Step 3: $f_\R$ is a smooth fiber bundle map:}\\
We  get as a direct consequence of Step 2 and Ehresmann's fibration theorem
that $f_\R$ is a smooth fiber bundle map.\\

\noindent {\bf Step 4: $f_\C$ is a smooth homotopy fibration:}\\
In view of Steps 1, 2 and 3, we may, and we will, assume that $f_\R$ has no critical
points. It remains to handle the atypical (including critical) values.
We assume that the atypical values
$z_1\, , \cdots\, ,z_k$ of $f_\C$
are isolated points in $\C^\ast$. Let $D_1\, , \cdots\, ,D_k$ be
small (disjoint) analytic neighborhoods of $z_1\, , \cdots \, ,z_k$ respectively.
We join $D_i$ by non-intersecting simple arcs $\alpha_i$ to $S^1$ and define
$$K\,:= \,S^1 \bigcup (\bigcup_{i=1}^k \alpha_i) \bigcup (\bigcup_{i=1}^k D_i)$$
and $W\,:=\, f^{-1}(K)$. If some $z_i \in S^1$, we assume that $\alpha_i$ is the constant arc at $z_i$.
Then $X_\C$ deformation retracts onto $W$, and hence $W$ has the same homotopy type as
$X_\R$ (since $X$ is a good complexification). 

The relative homotopy type of $(W\, , f^{-1}(S^1))$ is then given by the local 
topology changes at the critical fibers $f^{-1} (z_i)$ (cf. \cite{lam}). 
We shall now need to combine the above observation with the conclusion of Step 2.

Let $F_\R$ denote the
fiber of $f_\R \,:\, X_\R \, \longrightarrow\, S^1$, and let $F_\C$ denote the general
fiber of $f_\C$. Then $F_\C$ is also the general fiber of $\widehat f$. Pass to the cover $\til{X_\C}$
of $X_\C$ corresponding to $\pi_1 (F_\R)$ (so $\til{X_\C}$ is the fiber product of
$X_\C$ with the universal cover of ${\mathbb C}^*$), and let
$$P\,:\, \til{X_\C}\, \longrightarrow\, X_\C$$
denote the covering. Then $P^{-1} (X_\R)$ is homeomorphic to $F_\R \times \R$. We note that
$f_\C$ lifts to
a map $$\til{f_\C}\,:\, \til{X_\C}\, \longrightarrow\, \C\, ,$$ where $\C$ is
identified with the universal cover of $\C^\ast$,
and $\til{f_\C}$ restricts to a map $\til{f_\R}\,:\, \til{X_\R}\, \longrightarrow\, \R$,
where $\til{X_\R} \,=\,P^{-1} (X_\R)$, and $\R$ is identified with the
universal cover of $S^1$.

Fix fibers $F_\R \,\subset\, F_\C \,\subset\, \til{X_\C}$ as above. Then the quotient space
$\til{X_\C}/F_\R$ is contractible by the definition
of a good complexification. Further, $F_\C$ has the homotopy type of a finite $CW$
complex $H$ because it is a quasiprojective variety. 
Without loss of generality, we may assume that $H$ is a compact subset of $F_\C$.
Let $F$ denote the (compact finite) quotient complex $H/F_\R \,\subset\, \til{X_\C}/F_\R$. 
Each $z_i\,\in\, \C^\ast$ lifts to a countably infinite collection of atypical points
in $\C$. Hence (in its cell complex structure) $\til{X_\C}/F_\R$ has the same homotopy type as 
$H/F_\R$ with infinitely many cells attached, one for every lift of $z_i$, $i\,=\,
1,\cdots, k$. Since $\til{X_\C}/F_\R$ is contractible, it
follows that
\begin{enumerate}
 \item $k\,=\,0$, and
\item $H/F_R$ is contractible, in particular, $F_\C$ has the same homotopy type as $F_\R$.
\end{enumerate}

We summarize the conclusion of Step 4  as follows: \\
The map $f_\C$ is a smooth homotopy fibration over $\C^\ast$, where all fibers are
smooth manifolds of the same homotopy type. \\

\noindent {\bf Step 5: The fiber of $f_\R$ admits a good complexification:}\\
Let $N$ denote a general fiber of $f_\R$. The antiholomorphic involution $A$ preserves the
fibers $f_\C^{-1} (p)$ for $p \,\in\, S^1$. Let $X_p \,=\, f_\C^{-1} (p)$ and
$N_p \,=\, f_\R^{-1} (p)$. Then $N_p$ is diffeomorphic to $N$ for $p \in S^1$. It
suffices to show that $N_p$ is homotopy equivalent to $X_p$.

Now, by Step 4, the map $f_\C$ is a smooth homotopy fibration over
$\C^\ast$, where all fibers are smooth manifolds of the same homotopy type as
$X_p$. Hence we get the following commutative diagram of homotopy exact
sequences; the vertical arrows in the diagram are maps induced by inclusion:
$$
\begin{matrix}
\cdots \longrightarrow \, \pi_i({M}) & \longrightarrow & \pi_i({S^1}) &
\longrightarrow & \pi_{i-1}(N) & \longrightarrow & 
\pi_{i-1}(M) & {\longrightarrow} & \pi_{i-1}(S^1) & {\longrightarrow} \cdots\\
\,\,\,\,\quad\quad\,\,\,\,\,\Big\downarrow && \Big\downarrow && 
\Big\downarrow && \Big\downarrow && \Big\downarrow\\
\cdots \longrightarrow ~ \pi_i({X}) & \longrightarrow & \pi_i({\C^\ast}) & 
\longrightarrow &\pi_{i-1}(X_p) & \longrightarrow & 
\pi_{i-1}(X)&{\longrightarrow}& \pi_{i-1}({\C^\ast}) &{\longrightarrow} \cdots \\
\end{matrix}
$$
Since the homomorphism $\pi_i({S^1})\,\longrightarrow\, \pi_i({\C^\ast})$ is an
isomorphism for all $i$ (in fact they vanish for $i\, >\, 1$), and since the
homomorphism $\pi_i({M})\,\longrightarrow\, \pi_i({X})$ is an isomorphism by the
definition of a good complexification, it follows
that the inclusion of $N$ into $X_p$ induces isomorphisms $\pi_i({N})\,\longrightarrow\,
\pi_i({X_p})$ for all $i>1$ (by the 5-Lemma, using the fact that higher homotopy groups are abelian). Note further that the homomorphism
$\pi_1({N})\,\longrightarrow\, \pi_1({X_p})$ coincides with the identifications
of the kernels of $\pi_1(M)\,\longrightarrow\,\pi_1(S^1)$ and
$\pi_1(X)\,\longrightarrow\,\pi_1(\C^\ast)$. Also,
$\pi_0({N})\,\longrightarrow\, \pi_0({X_p})$ is an isomorphism.
Hence by the Whitehead Theorem, the inclusion of $N$ into $X_p$ is a homotopy
equivalence. This concludes Step 6 and proves the Theorem.
\end{proof}

We point out a couple of consequences of the above 
proof of Theorem \ref{cstar}. 

By Step 3, the map 
$ f_\R$  has no critical points, and by Step 6, $f$ has no critical 
points in $f^{-1}(S^1)$. Further, by Step 6, the fiber $F_\C$ 
has the same homotopy type as $F_\R$. Thus we have the following:

\begin{cor}\label{gcbdl0}
 Let $f\,:\, X \,\longrightarrow\, \C^\ast$ be a good complexification map between good
complexifications $X, \C^\ast$, and denote the antiholomorphic involution of $X$ by
$A_X$. Then
\begin{itemize}
\item $A_X(f^{-1}(p))\,=\, f^{-1}(p)$ for all
$p \,\in\, S^1$, and

\item{} $A_X(f^{-1}(p))$ is a good complexification of 
$f^{-1} (p)\bigcap X_\R$.
\end{itemize}
\end{cor}

\begin{proof} 
The only point to note is that $A_X$ preserves the fibers $f^{-1} (p)$. This in turn
follows from the condition $f \circ A_X\,=\, A_X \circ f$ in the
definition of a good complexification map.
\end{proof}

Combining Theorem \ref{cstar} and Corollary \ref{gcbdl0}, we obtain the following:

\begin{cor}\label{gcbdl}
Let $M$ be a closed smooth manifold admitting a good complexification. Then there exists a
good complexification map $f\,:\, X \,\longrightarrow\, \C^\ast$ between good
complexifications $X$ (of $M$) and $ \C^\ast$ (of $S^1$) such that $f$ is a smooth
homotopy fibration. Further, 
\begin{enumerate}
\item $M$ admits a smooth fibration over $ S^1$ with fiber $N$ such that $f_\R$ is a
projection onto the base $S^1$, and

\item $X$ admits a homotopy fibration over $ \C^\ast$ with generic fiber $Y$ such that
$f$ is a projection onto $\C^\ast$ and $Y$ is a good complexification of $N$.
\end{enumerate}
\end{cor}

\subsection{The Albanese}

We refer the reader to \cite{Se}, \cite{It}, \cite{SzSp} for
quasi-Albanese maps (see also \cite{Fu}), and to
\cite{nwy} for details on the semiabelian varieties (also called quasi-abelian varieties).
For a smooth quasiprojective variety $Z$, the quasi-Albanese map
\begin{equation}\label{az}
a_Z\, :\, Z\, \longrightarrow\, \text{Alb}(Z)
\end{equation}
is constructed by
choosing a point $z_0\, \in\, Z$ (the map on the degree zero part of ${\rm CH}_0(Z)$ does
not depend on the choice of the point $z_0$). We will replace the
quasi-Albanese $\text{Alb}(Z)$ by a torsor for it so that the construction of the
map from $Z$ does not require choosing a base point. Introduce the following equivalence
relation on $Z\times \text{Alb}(Z)$:
$$
(z_1\, , \alpha_1)\,\equiv\,
(z_2\, , \alpha_2) \ \ \text{ if }\ \ a_Z(z_1-z_2) \,=\, \alpha_2-\alpha_1\, ;
$$
as mentioned above, the map $a_X$ on the degree zero part of ${\rm CH}_0(Z)$ does
not depend on the choice of base point, so $a_Z(z_1-z_2)$ does not depend on $z_0$.
The corresponding quotient $\text{Alb}(Z)^1$
of $Z\times \text{Alb}(Z)$ is a torsor for $\text{Alb}(Z)$. Now consider the map
\begin{equation}\label{az2}
Z\, \longrightarrow\, \text{Alb}(Z)^1\, , \ x\, \longmapsto\, (z_0\, , a_Z(x))\, ,
\end{equation}
where $a_Z$ is the map in \eqref{az} constructed using $z_0$. It is straight-forward
check that this map does not depend on the choice of $z_0$. Henceforth, by
the quasi-Albanese of $Z$ we will mean $\text{Alb}(Z)^1$ constructed above, and
by quasi-Albanese map for $Z$ we will mean the map in \eqref{az2}.

If $Z$ is defined over $\mathbb R$, then both $\text{Alb}(Z)^1$ and $a_Z$ in
\eqref{az2} are also defined over $\mathbb R$.

\begin{lemma}
\label{alb}
Let $X_\C$ be a good complexification of $X_\R$. Then 
the quasi-Albanese of $X_\C$ is $(\C^\ast)^n$, where $n\,=\,b_1(X_\C)$.
\end{lemma}

\begin{proof}
Let $W_\C$ denote the quasi-Albanese of $X_\C$, and let
$$\al\,:\, X_\C \,\longrightarrow\, W$$ be the Albanese map. 
The variety $W_\C$ being defined over $\mathbb R$, is equipped with an
antiholomorphic involution $A_a$, and $\al$ intertwines the antiholomorphic involutions
of $X_\C$ and $W$, or in other words, we have
\begin{equation}\label{g-1}
A_a \circ \al\,=\ \al \circ A_X\, .
\end{equation}
Let $W_\R\,\subset\, W_\C$ denote the fixed-point locus of the
antiholomorphic involution. We have $\al(W_\C)\, \subset\, W_\R$.

Since the Albanese $W_c$ is a semi-abelian variety, there is an exact sequence $$1 
\,\longrightarrow\, (\C^\ast)^m \,\longrightarrow\, W_\C \,\longrightarrow\, Q 
\,\longrightarrow\, 1\, ,$$ where $Q$ is an abelian variety.

The antiholomorphic complex conjugation $A_a$ on $W_\C$ descends to an antiholomorphic
complex conjugation $A_q$ on $Q$ whose fixed
point set has dimension $\dim_\C (Q)$
(note that the fixed point set is nonempty because the identity element is fixed). Hence
$\al(X_{\R})$ has dimension at most $m + \dim_\C (Q)$ which is less than $\dim_\C (W_\C)$
unless $\dim_\C (Q)\,=\,0$. Since $\dim_\C (W_\C)\,=\,\dim_\R (W_\R)$, it follows that
$\dim_\C (Q)\,=\,0$. Hence the quasi-Albanese of $X_\C$ is $(\C^\ast)^n$, where
$n\,=\,b_1(X_\C)$.
\end{proof}

\begin{cor}
 \label{albgc} Let $X_\C$ be a good complexification of $X_\R$. Let $\al$ denote the Albanese map. Let $$\Pi\,:\,(\C^\ast)^n \,\longrightarrow\, \C^\ast$$
be a projection of the Albanese $W_\C\,=\,(\C^\ast)^n$ onto one of the factors. Then $\Pi \circ \al$ is a good complexification map.
\end{cor}

\begin{proof}
The quasi-Albanese of $X_\C$ is $(\C^\ast)^n$ by Lemma \ref{alb}, and
$\Pi$ is clearly a good complexification map. Therefore, from \eqref{g-1}
It follows that $\Pi \circ \al$ intertwines the antiholomorphic involution $A_X$ on $X_\C$
and the standard antiholomorphic involution $A_c$ on $\C^\ast$. Consequently, 
$\Pi \circ \al$ is a good complexification map.
\end{proof}

\subsection{Albanese Fibrations} The next lemma is a generalization of Case 3 B of
Proposition 4.3 of \cite{bm-imrn}.

\begin{lemma}\label{trivial1}
 Let $f\,:\, X \,\longrightarrow\, \C^\ast$ be a fibration that is 
a good complexification map between good complexifications $X, \C^\ast$ such that the general fiber
is $(\C^\ast)^m$. Then $f$ is a trivial fibration.
\end{lemma}

\begin{proof}
First note that $X$ may be regarded as a principal $(\C^\ast)^m$--bundle over $\C^\ast$.
Let $V$ denote the vector bundle on $\C^\ast$ of rank $m$ associated to the principal
$(\C^\ast)^m$--bundle $X$ for the natural action of $(\C^\ast)^m$
on ${\C}^{\oplus m}$. So $V$ splits as a
direct sum of $m$ line bundles $L_i$, i.e., $V \,=\, \bigoplus_{i=1}^m L_i$.
It suffices to show that $L_i$ is a trivial line bundle or equivalently
that the first Chern class $c_1(L_i)\,=\, 0$ for all $i$. 

Take any algebraic line bundle $L$ over $\C^\ast$.
The line bundle $L$ extends to an algebraic line bundle over the
${\mathbb P}^1$. To see this,
take the image in ${\mathbb P}^1$ of any divisor
in $\C^\ast$ representing $L$ by the homomorphism defined by
the inclusion map $\iota\, :\,
\C^\ast \, \hookrightarrow\, {\mathbb P}^1$.
Let $L'\,\longrightarrow\, {\mathbb P}^1$ be an extension of $L$. 
Therefore, $c_1(L)\,=\, \iota^* c_1(L')$. But
$$
\iota^*(H^2({\mathbb P}^1\, {\mathbb Z}))\,=\, 0\, .
$$
Therefore, we conclude that $c_1(L)\,=\, 0$.
\end{proof}

Essentially the same proof works when the base $\C^\ast$ is replaced by
$(\C^\ast)^n$. The only change is that we choose $({\mathbb P}^1)^n$ to be the
compactification of $(\C^\ast)^n$. This gives us the following:

\begin{cor}\label{trivial2}
Let $f\,:\, X \,\longrightarrow\, (\C^\ast)^n$ be a fibration that is 
a good complexification map between good complexifications $X, (\C^\ast)^n$ such that the general fiber
is $(\C^\ast)^m$. Then $f$ is a trivial fibration.
\end{cor}

We shall need the notion of a relative Albanese in the following. Let $f
\,:\, X\,\longrightarrow\, Y$ be a fibration. Then the relative Albanese
is a fibration $f_\alpha \,:\, X_\alpha\,\longrightarrow\, Y$, where the fiber $f_\alpha^{-1}(p)$ is the Albanese of 
$f^{-1}(p)$ for $p \,\in\, Y$. Note 
that the relative Albanese is a bundle of torsors, with canonically defined fibers, and hence in our case, a principal
$(\C^\ast)^k-$bundle for some integer $k\geq 0$ (multiplication is well-defined
in local trivializations and the fact that the fibers are canonically defined shows that multiplication agrees on overlaps). 
The following is the main technical tool of this paper.

\begin{theorem} \label{albbdl}
Let $M $ be a closed manifold admitting a good complexification.
Then there exists a good complexification $X_\C$ of $M\,=\, X_\R$, such that the
following holds: \\ Let $W_\C$ be the quasi-Albanese of $X_\C$, and 
let $\al\,:\, X_\C \,\longrightarrow\, W_\C$ be the quasi-Albanese map. Then $\al$ is a
homotopy smooth fibration with no critical points.
Further, if $F$ denote the general fiber, then
\begin{itemize}
\item $b_1(F)\,=\,0$, and

\item $F$ is a good complexification of the intersection $F\bigcap X_\R$.
\end{itemize}
\end{theorem}

\begin{proof}
If $W_\C$ is zero dimensional, then there is nothing to prove, so assume that $\dim 
W_\C\, \geq\, 1$.

We now proceed by induction on $b_1 (X_\C)$. By Lemma \ref{alb}, we have $W_\C\,=\,(\C^\ast)^n$.
Let $\Pi\,:\,(\C^\ast)^n\,\longrightarrow\,\C^\ast$ be a projection of the Albanese
$W_\C\,=\,(\C^\ast)^n$ onto one of the factors. Then $\Pi\circ \al$ is a good complexification map
by Corollary \ref{albgc}. 

Hence, by Theorem \ref{cstar}, the composition is a homotopy smooth fibration with no critical points.
Further, by Corollary \ref{gcbdl}, the fiber of $\Pi \circ \al$ over any 
point of $S^1$ is a good complexifications of the intersection of the fiber with $X_{\mathbb R}$. 

Then $\Pi \circ \al\,:\,X_\C\,\longrightarrow\,\C^\ast$ induces a good complexification map with fiber $F_z$ (say) over $z$.
Let $V_z$ denote the Albanese of the fiber $F_z\,=\,(\Pi \circ \al)^{-1}(z)$, and let
$$\widehat{\Pi \circ \al}\, :\, V\,\longrightarrow\, \C^\ast$$ be the bundle
$\{ V_z\,\mid\, z \in \C^\ast \}$.
For any such fibration, we have $b_1(V_z)\,\geq\,( b_1(X_\C) -1 )$. 

By Theorem \ref{cstar}, the above map
 $\widehat{\Pi \circ \al}$ defines a homotopy smooth
 fibration and the fibers $V_z$ are of the form $(\C^\ast)^k$ by Lemma \ref{alb}, where
$k\,=\,b_1(V_z)$. Since $V$ is a trivial bundle by Lemma \ref{trivial1}, we have 
$(k+1)\,\leq\,b_1(X_\C)$. It follows that $k\,=\,( b_1(X_\C) -1 )$.

By induction on the number of factors of $\C^\ast$ it follows that $F_z\,\longrightarrow\, V_z$ is a 
homotopy smooth fibration with no critical points. Hence $\al$ is a 
homotopy smooth fibration with no critical points.

Let $F$ denote the general fiber of $\al$.
It remains to show that $b_1(F)\,=\,0$. As above, pass to the relative Albanese bundle $W_F$ of Albanese tori. Then $W_F$
is a trivial smooth fibration over $W_\C\,=\,(\C^\ast)^n$ by Lemma \ref{trivial2} and hence has $b_1(W_F)\,>
\,n$ unless $b_1(F)\,=\,0$. Further, since $W_F$
is a semi-abelian variety, we have $b_1(W_F)\,\leq\,n$. Therefore, it now follows that $b_1(F)\,=\,0$.

That $F$ is itself a good complexification
now follows from Corollary \ref{gcbdl} by induction the number of factors of $\C^\ast$.
\end{proof}

The virtual first Betti number $vb_1(M)$ is defined to be the supremum of the first Betti numbers of finite sheeted
covers of $M$. It is known, \cite[Lemma 4.1]{bm-imrn}, that if $M$ admits a good complexification and $M_1$ is a finite-sheeted
cover of $M$, then $M_1$ admits a good complexification. We have the following immediate Corollary of Theorem \ref{albbdl}:

\begin{cor}
 \label{vb1} Suppose $M$ is a closed $n$--dimensional manifold admitting a good complexification. Then $vb_1(M)\,\leq\,n$.
\end{cor}

\begin{proof}
Let $M_1$ be a finite-sheeted cover of $M$. Then $M_1$ admits a good complexification \cite[Lemma 4.1]{bm-imrn}. Let $X$ be
a good complexification of $M_1$. By Theorem \ref{albbdl} this $X$ admits a smooth homotopy fibration with no critical points over
the Albanese $W$ of $X$. In particular, we have
$$
n\,=\, \dim_\C (X) \geq b_1(X)\,=\,b_1(M_1)\, .
$$
Since $M_1$ is arbitrary, the corollary follows.
\end{proof}

\section{Good Complexification Groups}\label{gcgrest}

\begin{theorem} \label{gcg}
Let $G$ be a good complexification group. Then there exists a finite index subgroup $G_1$ of $G$ such that 

\begin{enumerate}
 \item There is an exact sequence: $$1\,\longrightarrow\, H\,\longrightarrow\, G_1
\,\longrightarrow\, \Z^k \,\longrightarrow\, 1\, ,$$ where $k$ can be zero.
\item $H$ is a finitely presented good complexification group with $vb_1(H)\,=\,0$, where
$vb_1(H)$ denotes the virtual first Betti number of $H$.
\end{enumerate}
\end{theorem}

\begin{proof}
Let $M$ be a closed smooth $n$--manifold admitting a good complexification with $\pi_1(M)\,=\,G$. By Corollary \ref{vb1}, 
there is a finite sheeted cover $M_1$ of $M$ such that
$$vb_1(M)\,=\,b_1(M_1)\,=\,m\,\leq\,n\, .$$ Also by Theorem \ref{albbdl}, $M_1$ fibers over $(S^1)^m$ with fiber
a closed $(n-m)$--manifold $N$. Further $N$ admits a good complexification by 
Theorem \ref{albbdl}. Hence $vb_1(N)\,=\,0$; indeed, if
$vb_1(N)\,\not=\,0$, by applying Theorem \ref{albbdl} again, $vb_1(M)\,=\,vb_1(M_1)\,>\,m$.
Defining $H\,=\,\pi_1(N)$ and $G_1 \,=\,\pi_1(M_1)$,
 we obtain the following homotopy long exact sequence: $$1 \,\longrightarrow\, H
 G_1 \,\longrightarrow\, \Z^{k'}\, .$$
Hence the theorem follows.
\end{proof}

\begin{cor} \label{not}
 The following are not good complexification groups:
\begin{enumerate}

\item Groups with infinite $vb_1$, in particular, large groups.

\item Hyperbolic {\rm CAT(0)} cubulated groups; in particular,
infinite hyperbolic Coxeter groups.

\item Solvable groups that are not virtually abelian.

\item 2-manifold groups that are not virtually abelian.

\item 3-manifold groups that are not virtually abelian.
\end{enumerate}
\end{cor}

\begin{proof} (1)~Groups with infinite $vb_1$ cannot be good complexification groups by Theorem \ref{gcg}. Recall that
$G$ is large if $G$ has a finite index subgroup that surjects onto the free group $F_2$ of
two generators \cite{bm-imrn}. Clearly, large groups have infinite $vb_1$.\\ 

(2) We refer the reader to \cite{gr-hg, farb, gr-ai} for generalities on hyperbolic and relatively hyperbolic groups,
to \cite{sageev, wise} for generalities on CAT(0) cube-complexes and their properties and to \cite{hw} for the theory of virtually special
complexes.

A hyperbolic CAT(0) cubulated group $G$ is virtually special and hence has infinite
$vb_1(G)$ \cite{agol, wise}. 
We are thus reduced to Case (1) above.

The second statement follows from \cite{hw} where the
authors show that infinite Coxeter groups are CAT(0) cubulated virtually special.\\

Non-elementary hyperbolic, or relatively hyperbolic, groups have asymptotic cones with cut points. 
On the other hand, a good complexification group with positive $vb_1$ must be virtually of the form $G_1\,=\,\Z^k \times H$
by Theorem \ref{gcg}. Such a group has as asymptotic cone
$$AC(G)\,=\,AC(G_1)\,=\,\R^k \times AC(H)\, ,$$ where $AC(H)$ is an asymptotic cone of $H$, and $k\,\geq\,1$.
Hence $AC(G)$ can have a cut-point if and only if $H$ is finite and $k\,=\,1$, which means that $G$ is virtually cyclic.\\

(3) A solvable group $G$ that is not virtually abelian admits, up to finite index, an exact sequence $$1\,\longrightarrow\, N\,\longrightarrow\, G
\,\longrightarrow\, \Z^k\,\longrightarrow\, 1\, ,$$
where $k$ is maximal, and $b_1(N)\,>\,0$. This contradicts Theorem \ref{gcg}. Note that we avoided using \cite{an} in this context. \\

(4) Since 2-manifold groups $G$ that are not virtually abelian have $vb_1(G)\,=\,\infty$, 
this follows from Case 1 above.\\

(5) By work of a large number of people culminating in the work of Agol 
and Wise \cite{agol, wise} (see \cite{afw}, especially Diagram 1 in page 36 for details), a 3-manifold group $G\,=\,\pi_1(M^3)$ is large 
 unless $G$ is virtually solvable,
or equivalently, $M^3$ has geometry modelled on one of $S^3$, $E^3$, $S^2 \times \R$, $Nil$ and $Sol$. Now, by Case 2 above, the group $G$ must be
virtually abelian. Hence $M$ must have geometry modelled on one of $S^3$, $E^3$ 
and $S^2\times \R$. 
\end{proof}

\begin{rmk}\mbox{}
\begin{enumerate}
\item Case 3 of Corollary \ref{not} implies that any unipotent representation of a good 
complexification group is virtually abelian.

\item It also follows from Theorem \ref{gcg} that the commutator subgroup of a finite index subgroup of a
good complexification group is finitely presented.

\item Several cases of Corollary \ref{not} above may be strengthened by saying that the groups listed there cannot 
even appear as quotients of good complexification groups. The proof involves applying Theorem \ref{gcg} as done in 
Corollary \ref{not}.
\end{enumerate}
\end{rmk}

\section{Examples and Applications} \label{ld} In this Section, we apply Theorem \ref{albbdl} to furnish conclusions for low-dimensional
examples.

\subsection{2 Manifolds}\label{2m}

\begin{lemma} Question \ref{tqn} has an affirmative answer for 2-manifolds. \label{2mfd}
 \end{lemma}

\begin{proof}
 By Case 4 of Corollary \ref{not}, the only 2-manifolds that can possibly admit good complexifications are those finite covered by
the sphere
$S^2$ or the torus $S^1 \times S^1$. That all such manifolds do admit metrics of non-negative curvature is classical.
Hence Question \ref{tqn} has an affirmative answer for 2-manifolds.
\end{proof}

\subsection{3 Manifolds}
\label{3m}

\begin{theorem}[\cite{bm-imrn}] Question \ref{tqn} has an affirmative answer for 3-manifolds. \label{3mfd}
 \end{theorem}

\begin{proof}
 By Case 5 of Corollary \ref{not}, the only 3-manifolds that can possibly admit good complexifications are those finite covered by
$S^3, S^2 \times S^1$ or $S^1 \times S^1 \times S^1$. That all such manifolds do admit metrics of non-negative curvature is 
a consequence of the Geometrization Theorem of Perelman.
\end{proof}

\subsection{Simply Connected 4-manifolds} \label{totcom} A stronger notion of good
complexification appears implicitly in the work of McLean \cite{mmc}. There, the author
observes first that any smooth complex affine variety $X$ has a natural structure of a
symplectic manifold (up to symplectomorphisms) when regarded as a real smooth manifold. The
cotangent bundle $T^\ast M$ of a closed manifold $M$ carries a natural symplectic structure
(it is called the Liouville symplectic form). Corollary 1.3 of \cite{mmc} implies that if
$X$ is {\em symplectomorphic} to $T^\ast M$ for some closed manifold $M$ (in particular, $M$
admits a good complexification in the sense of \cite{tot} (recalled in Section \ref{se1})), 
then $M$ is rationally elliptic. We say that $M$ admits a {\bf symplectic good complexification} if $T^\ast M$ with its 
natural symplectic structure is symplectomorphic to a smooth complex affine variety $X$ with its
natural underlying symplectic structure.

\begin{prop}\label{sc4}
 Let $M$ be a closed simply connected 4-manifold admitting a symplectic good complexification. Then $M$ admits
a metric of non-negative curvature (i.e., Question \ref{tqn} has an affirmative answer).
\end{prop}

\begin{proof}
 From \cite[Corollary 1.3]{mmc}, the manifold $M$ is rationally elliptic. Hence $b_2(M)\,\leq\, 2$. By the classification
of simply connected 4-manifolds, $M$ is homeomorphic to one of the following: \\ $$S^4,~ \C P^2,~ \C P^2\# \C P^2,
~\C P^2\# \bbar{ \C P^2}, ~ S^2 \times S^2\, .$$
That each of these admits
a metric of non-negative curvature follows from the examples constructed by Totaro in \cite{tot}.
\end{proof}

\subsection{4 manifolds with fundamental group $\Z$ and $\Z\oplus \Z$}\label{somnath}\mbox{}\\
We shall need the following theorem.

\begin{theorem}[Smale, Hatcher]\label{smalehatcher}\mbox{}
\begin{enumerate}
\item ${\rm Diff}(S^2)$ is homotopy equivalent to ${\rm SO}(3)$ \cite{sm}.

\item ${\rm Diff}(S^3)$ is homotopy equivalent to ${\rm SO}(4)$ \cite{ha}.
\end{enumerate}
\end{theorem}

\begin{cor} \mbox{} \label{shcor}
\begin{enumerate}
\item[(a)] An $S^2$ bundle over $T^2$ is trivial after (at most) passing to
a double cover of $T^2$.

\item[(b)] An $S^2$ bundle over $T^3$ is trivial after (at most) passing to
a finite-sheeted cover of $T^3$.

\item[(c)] An $S^3$ bundle over $T^2$ is trivial after (at most) passing to
a double cover of $T^2$.

\item[(d)] An $S^3$ bundle over $T^3$ is trivial after (at most) passing to
a finite-sheeted cover of $T^3$.
\end{enumerate}
\end{cor}

\begin{proof} \mbox{} \\
(a) The collection of distinct $S^2$-bundles over $T^2$ is in bijective correspondence with homotopy classes of maps $[T^2, {\rm BDiff}(S^2)]$. 
By the first statement of Theorem \ref{smalehatcher}
this is equivalent to $[T^2, {\rm BSO}(3)]$. By the homotopy long exact
sequence, $\pi_n({\rm BSO}(3))\,= \,\pi_{n-1} (SO(3))$. Since ${\rm BSO}(3)$ is simply connected, 
it follows that (after equipping $T^2$ with the standard CW complex structure consisting of one 0-cell, two 1-cells and one 2-cell) 
any map from $T^2$ to ${\rm BSO}(3)$ induces a map from $S^2$ to ${\rm BSO}(3)$, where $S^2$ is the quotient of $T^2$ obtained by collapsing its
one skeleton to a point. Hence $[T^2, {\rm BDiff}(S^2)] \,=\, \pi_2 ({\rm BSO}(3))\,= \,\pi_1 (SO(3)) \,=\, \Z/2\Z$. It follows that after passing to
a double cover of $T^2$ if necessary, the pullback $S^2$ bundle becomes trivial.\\

\noindent (b) Repeating the argument for (a) above, we find that any $S^2$ bundle over $T^3$ can be made trivial over the 2-skeleton
after passing to a 
finite-sheeted cover. Hence any map from $T^3$ to ${\rm BSO}(3)$ induces a map from $S^3$ to ${\rm BSO}(3)$, where $S^3$ is the quotient of $T^3$ obtained 
by collapsing its
two skeleton to a point. Hence $[T^3, {\rm BDiff}(S^2)]\,=\, \pi_3 ({\rm BSO}(3))
\,=\, \pi_2 ({\rm SO}(3)) \,=\, 0$.
It follows that after passing to
a finite-sheeted cover of $T^3$ if necessary, the pullback $S^2$ bundle becomes trivial.\\

\noindent (c)
We repeat the argument in (a), replacing ${\rm BSO}(3)$ with ${\rm BSO}(4)$ (and using the second statement of Theorem \ref{smalehatcher})
to obtain $[T^2, {\rm BDiff}(S^3)]\,=\, \pi_2 ({\rm BSO}(4))\,=\, \pi_1 ({\rm SO}(4))
\,=\, \Z/2\Z$. Here the last statement follows from the fact that $SO(4)$ is
diffeomorphic to $SO(3) \times S^3$. The rest of the argument is identical to that in (a).\\

\noindent (d) We repeat the argument in (b) to obtain $[T^3, {\rm BDiff}(S^3)]\,=\,
\pi_3 ({\rm BSO}(4))\,= \,\pi_2 ({\rm SO}(4)) \,=\, 0.$
The Corollary follows.
\end{proof}

\begin{prop}\label{c4}
 Let $M$ be a closed 4-manifold admitting a good complexification. Further suppose that $\pi_1(M)$ is infinite, torsion-free
abelian. Then $\pi_1(M)$ is isomorphic to 
 $\Z^d$, where $d=1, 2$ or $4$. Further, the manifold $M$ admits
a finite-sheeted cover with a metric of non-negative curvature (i.e., Question \ref{tqn} has an affirmative answer
up to finite-sheeted covering).
\end{prop}

\begin{proof} By Theorem \ref{albbdl}, the 4-manifold $M$ admits a fiber bundle structure over $(S^1)^d$, where $1 \leq d\leq 4$
(by Corollary \ref{vb1}). Let $N$ denote the fiber. Then $b_1(N)=0$ 
(by Theorem \ref{albbdl} again). Hence $N$ cannot be one-dimensional, and hence $d\,=
\,1, 2$ or $4$.

If $d\,=\,4$, then $M$ is diffeomorphic to $(S^1)^4$.

If $d\,=\,2$, then $M$ fibers over $(S^1)^2$ with fiber $N$ a 2-manifold admitting a good complexification and having $b_1(N)=0$. 
Hence $M$ is an $S^2$-bundle over $T^2$ and is finitely covered by $S^2 \times T^2$. Clearly, $S^2 \times T^2$ admits
a metric of 
non-negative curvature.

If $d=1$, then $M$ fibers over $S^1$ with fiber $N$ a 3-manifold admitting a good complexification and having $b_1(N)=0$. 
 By Theorem \ref{3mfd}, $N$ must be finitely covered by $(S^1)^3$ or $S^1 \times S^2$ or $S^3$. Hence after passing to
a finite-sheeted cover if necessary, $M$ is a fiber bundle over $B$ with fiber $N_1$, where $B, N_1$ satisfy one
of the following (Corollary \ref{shcor}):

\begin{enumerate}
 \item $B\,=\, (S^1)^4$, and $N_1$ is a point. In this case, $M$ admits a finite-sheeted cover diffeomorphic to $(S^1)^4$.

\item $B\,= \,(S^1)^2$, and $N_1 \,=\, S^2$. In this case, $M$ admits a finite-sheeted cover diffeomorphic to $(S^1)^2 \times S^2$
as for the case $d=2$ above.

\item $B\,=\, S^1$, and $N_1 \,=\, S^3$. In this case, $M$ admits a finite-sheeted cover diffeomorphic to $S^1 \times S^3$.
\end{enumerate}
The Theorem follows.
\end{proof}

\subsection{5 Manifolds} The argument for Proposition \ref{c4} coupled with Corollary 
\ref{shcor} furnishes the following fact for 5-manifolds as well.

\begin{prop}\label{c5}
 Let $M$ be a closed 5-manifold admitting a good complexification. Further suppose that $\pi_1(M)$ is infinite, torsion-free
abelian. Then $\pi_1(M)$ is isomorphic to 
 $\Z^d$, where $d=1, 2, 3$ or $5$. Further, if $d \,=\, 2, 3$ or $5$, then $M$ admits
a finite-sheeted cover $M_1$ homeomorphic to $S^2 \times T^3$ or $S^3 \times T^2$ or $T^5$.
In particular, $M_1$ admits a metric of non-negative curvature.
\end{prop}

We do not repeat the proof of Proposition \ref{c4} here but observe that (after passing to a finite-sheeted cover),
$d\,=\,2$ corresponds to an $S^3$-bundle over $T^2$, $d=3$ corresponds to an $S^2$-bundle over $T^3$,
$d\,=\,5$ corresponds to $T^5$.

\subsection{Branched Covers}

\begin{prop} Let $f\,:\, X \,\longrightarrow\, Y$ be a good complexification map between good complexifications
$X, Y$ such that $f$ is a branched cover. Then $f$ is actually \'etale.
\label{etale} \end{prop}

\begin{proof} The branch locus of $f_\C$ has to be of complex codimension one in $Y_\C$. Hence the branch locus 
$B_\R\,\subset\,Y_\R$ of $f_\R$ is of real codimension one. If $B_\R\,\neq\,\emptyset$ (or equivalently, if
$f_\R$ is not \'etale) then $X_\R$ cannot be a manifold because
the local model of $X_\R$ around a preimage of a point in $B_\R$
will consist of more than two half-spaces glued along a codimension one hyperplane. It follows that $B_\R\,=\,\emptyset$
and hence $f_\R$ is \'etale.

It remains to show that $f_\C$ is \'etale. First, since $f_\R$ is \'etale, and $Y_\R\,\subset\, Y_\C$ is a homotopy
equivalence, we can construct a new good complexification $X_\C^1$ of $X_\R$ by simply taking $X_\C^1$ to be
the cover of $Y_\C$ corresponding to $f_{\R, \ast}(\pi_1(X_\R))$. Hence by lifting the map $f_\C\,:\, X_\C\,
\longrightarrow\, Y_C$ to $X_\C^1$, we have an algebraic homotopy equivalence $h$
between $X_\C$ and $X_\C^1$ which is the identity on $X_\R$. Since $X_\R$ is Zariski-dense in $X_\C$, $h$ must be
the identity map. i.e., $X_\C\,=\,X_\C^1$.
\end{proof}

\subsection{Questions}
In order to answer Question \ref{tqn} affirmatively, it suffices by Theorems \ref{albbdl}, \ref{gcg}, to answer the following questions affirmatively: 

\begin{qn}\label{tqn1} \mbox{}
\begin{enumerate}
\item[(a)]Suppose that the virtual first Betti number of a closed manifold $M$ is zero: $vb_1(M)\,=\,0$.
If $M$ admits a good complexification, then does $M$ admit a metric of non-negative curvature?
\item[(b)] If $M$ is a fiber bundle over $(S^1)^d$ admitting a good complexification, is the bundle trivial (at least up to finite sheeted covers)?
\end{enumerate}
\end{qn}

We divide Question \ref{tqn1} (a) into two further questions: 

\begin{qn} \label{tqn2}\mbox{}
\begin{enumerate}
\item[(a)] Suppose that $G$ is a good complexification group with $vb_1(G)\,=\,0$, is $G$ finite?

\item[(b)] If a closed manifold $M$ with finite fundamental group admits a good complexification, then does $M$
admit a metric, at least virtually, of non-negative curvature?
\end{enumerate}
\end{qn}

\begin{rmk} An affirmative answer to both parts of Questions \ref{tqn1} and \ref{tqn2} exist if and only if Question \ref{tqn} has 
an affirmative answer.\label{eqv} \\
To see this, note that if $M$ admits a metric of non-negative curvature, then 
the universal cover splits (by the Cheeger-Gromoll splitting Theorem \cite[Theorem 69 , page 288]{pet}) as a metric 
product $N \times \R^k$ and hence $G=\pi_1(M)$ is virtually $\Z^k$. 

Conversely, let $G\,=\, \pi_1(M)$ such that $M$ 
admits a good complexification. By Theorem \ref{gcg}, $G$ is virtually of the form $H \rtimes \Z^k$, where $H\,=\, 
\pi_1(N)$ where $N$ admits a good complexification and $vb_1(N)\,=\,0$. By a positive answer to Question \ref{tqn2}(a), 
$H$ must be finite. Hence $N$ admits a finite simply connected cover $N_1$, which in turn admits a good 
complexification \cite[Lemma 4.1]{bm-imrn}. Hence $N_1$ admits a metric of non-negative curvature by a positive 
answer to Question \ref{tqn2}(b). Hence $N_1 \times (S^1)^k$ admits a metric of non-negative curvature. 
Consequently, by a positive 
answer to Question \ref{tqn1}(b), $M$ admits a finite sheeted cover with a metric of non-negative curvature.

Since every finite group is a subgroup of some ${\rm SU}(n)$, every finite group is a good complexification group
\cite{tot}, and also the fundamental group of a closed manifold of non-negative curvature.
\label{harish}\end{rmk}

\section*{Acknowledgments} The authors thank Burt Totaro, Harish Seshadri and Somnath Basu for 
useful email correspondence and conversations relevant to Section \ref{totcom}, Remark \ref{harish}
and Section \ref{somnath} respectively. Parts of this work were accomplished during visits of one
or more of the authors to Tata Institute of Fundamental Research (TIFR), Mumbai; Harish Chandra
Research Institute (HRI), Allahabad; RKM Vivekananda University, Belur; Indian Statistical
Institute (ISI), Kolkata. A substantial amount of the work was done during a Discussion Meeting on
Symplectic and Contact Topology held at HRI and TIFR in December 2014. We thank all these
institutes for their hospitality.

\end{document}